\documentclass{article}
\usepackage[utf8]{inputenc}
\usepackage[a4paper]{geometry}
\usepackage{hyperref,mathrsfs,enumitem,amsmath, amssymb, amsthm, graphicx,color,comment,array}
\usepackage{tikz-cd} 
\usepackage{setspace}
\numberwithin{equation}{section}
\usepackage{parskip}
\setlist{  
  listparindent=\parindent,
  parsep=0pt,
}

\DeclareMathOperator{\orb}{orb}

\newcommand{\N}{\mathbb{N}}

\newcommand{\C}{\mathbb{C}}
\newcommand{\f}{\frac 12}

\usepackage{mathtools}

\newtheorem{theorem}{Theorem}[section]

\theoremstyle{definition}\newtheorem{remark}[theorem]{Remark}
\theoremstyle{definition}\newtheorem{fact}[theorem]{Fact}
\theoremstyle{definition}\newtheorem{diagram}[theorem]{Diagram}
\theoremstyle{definition}\newtheorem{definition}[theorem]{Definition}
\theoremstyle{definition}\newtheorem{question}[theorem]{Question}

\theoremstyle{definition}\newtheorem{example}[theorem]{Example}

\begin{document}
\title{Strong topological transitivity, hypermixing, and their relationships with other dynamical properties}
\author{Ian Curtis\footnote{St. Olaf College, Northfield, MN 55057}, Sean Griswold\textsuperscript{$*$}, Abigail Halverson\textsuperscript{$*$},  Eric Stilwell\textsuperscript{$*$}, \\ Sarah Teske\textsuperscript{$*$}, David Walmsley\footnote{Department of Mathematics, Statistics and Computer Science, St. Olaf College, Northfield, MN 55057}, Shaozhe Wang\textsuperscript{$*$}}




\date{}
\maketitle

\renewcommand{\thefootnote}{\fnsymbol{footnote}} 
\renewcommand{\thefootnote}{\arabic{footnote}} 
\begin{abstract}
Recently, two stronger versions of dynamical properties have been introduced and investigated: strong topological transitivity, which is a stronger version of the topological transitivity property, and hypermixing, which is a stronger version of the mixing property. We continue the investigation of these notions with two main results. First, we show there are dynamical systems which are strongly topologically transitive but not weakly mixing. We then show that on $\ell^p$ or $c_0$, there is a weighted backward shift which is strongly topologically transitive but not mixing.
\end{abstract}


\section{Introduction}\label{Introduction}
A continuous linear map (henceforth an \textit{operator}) $T:X\to X$ on a topological vector space $X$ is called \textit{hypercyclic} if there exists an element $x\in X$ (called a \textit{hypercyclic element}) whose \textit{orbit under $T$}, given by $\orb(T,x)=\{x,Tx,T^2x,T^3x,\ldots\}$, is dense in $X$. Denote by $\N$ the set of positive integers and $\N_0$ the set of nonnegative integers. Given nonempty open subsets $U$ and $V$ of $X$, the \textit{return set} $N(U,V)$ is defined as $N(U,V)=\{n\in \N_0: T^n(U)\cap V\not =\emptyset\}$. If $X$ is a metric space, we say a continuous mapping $T:X\to X$ is \textit{topologically transitive} (resp. \textit{weakly mixing}; resp. \textit{mixing}) if for all nonempty open subsets $U$ and $V$ of $X$, $N(U,V)$ is nonempty (resp. contains arbitrarily long intervals; resp. is cofinite). When $X$ is second countable and Baire, hypercyclicity is equivalent to topological transitivity by the Birkhoff transitivity theorem \cite{Birkhoff1}. Hypercyclicity has become a very active area of research in operator theory that has many connections to other branches of mathematics. For more background on its history, central ideas, and a detailed introduction, we refer the interested reader to the monographs \cite{BayartMatheron,GrossePeris}. 

Given our definitions, the following implications are immediate:
\begin{center}
    \begin{tabular}{ccccc}
    mixing & $\implies$ & weakly mixing & $\implies$ & topologically transitive. 
\end{tabular}
\end{center}

\noindent Recently in \cite{Ansari1,Ansari2,Ansari3,Ansari4}, stronger versions of topological transitivity and mixing were introduced and studied.
\begin{definition}
Let $X$ be a topological vector space and $T$ be an operator on $X$. We say $T$ is \textit{strongly topologically transitive} (or \textit{strongly hypercyclic} when $X$ is second countable and Baire) if, for each nonempty open subset $U$ of $X$,
\begin{align*}
    X\setminus \{0\}\subseteq \bigcup_{n=0}^\infty T^n(U).
\end{align*}
We say $T$ is \textit{supermixing} if,  for each nonempty open subset $U$ of $X$,
\begin{align*}
    X = \overline{\bigcup_{i=0}^\infty \bigcap_{n=i}^\infty T^n(U)}.
\end{align*}
We say $T$ is \textit{hypermixing} if,  for each nonempty open subset $U$ of $X$,
\begin{align*}
    X = \bigcup_{i=0}^\infty \bigcap_{n=i}^\infty T^n(U).
\end{align*}
\end{definition}

\begin{remark}
Because of \cite[Remark 2.5]{Ansari4}, the definition of hypermixing that appears above coincides with the original definition of hypermixing in \cite[Definition 1.1]{Ansari4}.
\end{remark}
The following implications follow quickly from the definitions:

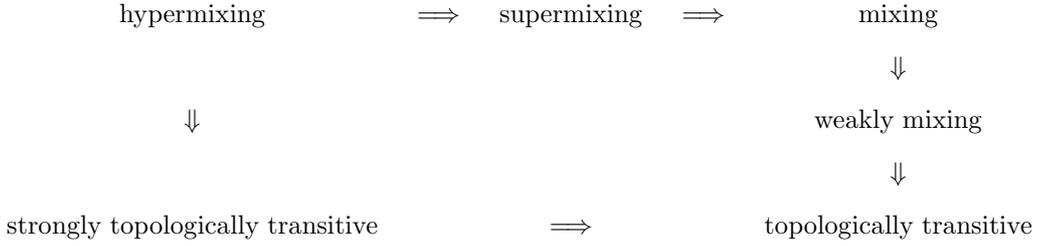
\begin{center}
\begin{tabular}{ccccc}
    hypermixing & $\implies$ & supermixing & $\implies$ & mixing\\
                &            &             &            & $\Downarrow$\\
    $\Downarrow$ &           &             &            & weakly mixing\\
                &            &             &            & $\Downarrow$\\
    strongly topologically transitive &  & $\implies$ && topologically transitive
\end{tabular}
\end{center}
In \cite{Ansari4}, the author explored several of the implications not yet present in the above diagram for the family of weighted backward shift operators on $\ell^p$ and $c_0$. For example, they showed that every mixing weighted backward shift on $\ell^p$ or $c_0$ is in fact supermixing, but not necessarily hypermixing. They also gave an example of a mixing weighted backward shift which is not strongly topologically transitive. 
Examples of weakly mixing weighted shifts which are not mixing are well-known \cite[Remark 4.10]{GrossePeris}, and it turns out that topological transitivity implies the weak-mixing property for the family of weighted backward shifts (see \cite[Section 4.1]{GrossePeris}). 
Given these results, there are two implications not present in the above diagram left to explore for the family of weighted backward shifts, which led the author in \cite{Ansari4} to pose the following questions.
\begin{question}
Let $B_w$ be a strongly topologically transitive weighted backward shift on $\ell^p$ or $c_0$. Must $B_w$ be hypermixing? Must $B_w$ be mixing?
\end{question}

In Section \ref{backward shifts}, we answer these questions in the negative with the following result.
\begin{theorem}\label{shift theorem}
    On $\ell^p$ or $c_0$, there exists a weighted backward shift which is strongly topologically transitive (i.e. strongly hypercyclic) but not mixing.
\end{theorem}
\noindent This result, along with those in \cite{Ansari4}, complete the above diagram for the family of weighted backward shift operators on $\ell^p$ or $c_0$.

Before getting to this, in Section \ref{general stuff}, we extend the definitions of strong topological transitivity, supermixing and hypermixing beyond the setting of operators on vector spaces to the setting of continuous maps on metric spaces. We then classify some classic dynamical systems as hypermixing or strongly topologically transitive. In doing so, we provide, perhaps surprisingly, an example of a strongly topologically transitive map which is not weakly mixing. We also show that there are mixing mappings which are not supermixing, which shows that the above diagram for metric spaces is complete as shown. 

\section{The metric space setting}\label{general stuff}

In this section, we take $X$ to be a metric space with metric $d$ and $T:X\to X$ a continuous mapping. The pair $(X,T)$ is called a \textit{dynamical system}, although we will simply call $T:X\to X$ a dynamical system. We adopt the convention in operator theory to write $Tx$ for $T(x)$.

The definitions of supermixing and hypermixing work perfectly well in this setting without any changes. To quickly find examples of hypermixing mappings, there is a simple criterion to check that covers several important dynamical systems.

\begin{theorem}\label{simple criterion}
    Suppose $T:X\to X$ is a dynamical system with the property that, for every nonempty open set $U$, there exists $j\in \N$ such that $T^j(U)=X$. Then $T$ is hypermixing.
\end{theorem}

\begin{proof}
    We first show that $T$ is surjective. Let $y\in X$ and $U$ be a nonempty open subset of $X$. Then $T^j(U)=X$ for some $j\in \N$, which implies $T^j x = T(T^{j-1}x)=y$. Hence $T$ is surjective, which implies $T^n(X)=X$ for all $n\in\N$. Hence for $n\geq j$, we have $T^n(U)=T^{n-j}(X)=X$, which implies
    \begin{align*}
    X=\bigcap_{n=j}^\infty T^n(U) \supseteq \bigcup_{i=0}^\infty \bigcap_{n=i}^\infty T^n(U),
\end{align*}
   showing that $T$ is hypermixing.
\end{proof}

\begin{example}
     \begin{enumerate}[label=\upshape(\alph*)]
        \item The ``tent map" $T:[0,1]\to [0,1]$ defined by $T(x)=2x$ if $x\in[0,\f]$, $T(x)=2-2x$ if $x\in ]\f,1]$ is easily seen to satisfy the conditions in Theorem \ref{simple criterion}; see \cite[Example 1.12(a)]{GrossePeris}.
        \item On the unit circle $S^1$, the ``doubling map" $T:S^1\to S^1$ given by $T(z)=z^2$ satisfies the above criterion; see \cite[Example 1.12(b)]{GrossePeris}.
        \item On the space of all $0$-$1$ sequences $\Sigma_2=\{(x_n)_{n\geq 0}: x_n \in \{0,1\}\}$, consider the ``shift" $\sigma: \Sigma_2\to \Sigma_2$ defined by
        \begin{align*}
            \sigma(x_0,x_1,x_2,\ldots)=(x_1,x_2,x_3,\ldots).
        \end{align*}
        The metric on $\Sigma_2$ is given by
        \begin{align*}
            d(x,y)=\sum_{n=0}^\infty \frac{|x_n-y_n|}{2^n}.
        \end{align*}
        It is easy to see that if $x=(x_n)_{n\geq 0}$ and $y=(y_n)_{n\geq 0}$ are sequences in $\Sigma_2$ with $x_n=y_n$ for $n=0,1,\ldots,j$, then $d(x,y)\leq 2^{-j}$. To show $\sigma$ has the property above, let $U$ be a nonempty open subset of $X$. Then there exists $u\in U$ and $\epsilon>0$ for which the ball $B_\epsilon(u)\coloneqq\{x\in X: d(x,u)<\epsilon\}\subseteq U$. 
        
        \sloppy Choose $j\in\N$ such that $2^{-j}<\epsilon$, and let $x\in X$. Then the element $v\coloneqq (u_0,\ldots,u_j,x_0,x_1,x_2,\ldots)$ belongs to $U$ since $d(v,u)\leq 2^{-j}<\epsilon$, and we have arranged $\sigma^j(v)=x$. Hence $\sigma^j(U)=X$, so $\sigma$ is hypermixing.
        \item No (linear and bounded) operator $T$ acting on a Banach space $X$ can satisfy the property in Theorem \ref{simple criterion}. Hence that criterion does not apply to a linear setting. 
     \end{enumerate}
\end{example}

It was shown in \cite[Remark 2.5]{Ansari4} that an injective operator on a topological vector space is never hypermixing. The argument given in \cite[Theorem 2.3]{Ansari4} was based on a hypermixing criterion that exploits the underlying linearity of the operator. We show that, despite the absence of linearity, no injective mapping on a metric space is even supermixing. 

\begin{theorem}\label{injective}
    Let $X$ be a metric space and $T:X\to X$ be a continuous map which is not the identity. If $T$ is injective, then $T$ is not supermixing.
\end{theorem}

\begin{proof}
    We establish something stronger than not being supermixing by showing there exists a nonempty open subset $U$ of $X$ such that $T^n(U)\cap T^{n+1}(U)=\emptyset$ for each $n\in \N_0$. The main ingredient to do this is the following fact.
    \begin{fact}
        If $T:X\to X$ is not the identity, then there exists a nonempty open subset $U$ of $X$ for which $U\cap T(U)=\emptyset$.
    \end{fact}
    Granting this fact, we conclude as follows. Let $U$ be a nonempty open subset of $X$ given by the fact. Let $n\in \N_0$. Towards a contradiction, suppose $T^n(U) \cap T^{n+1}(U)\not =\emptyset$. Then $T^n u_1 =T^{n+1}u_2$ for some $u_1,u_2\in U$. Since $T^n$ is injective, we have $u_1=Tu_2$. But then $T(U)\cap U\not =\emptyset$, which is impossible. Thus for any $n\in \N_0$, we must have $T^n(U)\cap T^{n+1}(U)=\emptyset$, which shows each set $\bigcap_{n=i}^ \infty T^n(U)$ is empty.
\end{proof}

\begin{proof}[Proof of the fact.] Let $x\in X$ with $Tx\not =x$, and let $\epsilon=d(x,Tx)$. Since $T$ is continuous, there exists $\delta>0$ for which $d(x,y)<\delta$ implies $d(Tx,Ty)<\frac{\epsilon}{2}$. Let $\rho=\min(\delta,\frac{\epsilon}{2})$ and $U$ be the ball $B_\rho(x)\coloneqq \{y\in X: d(x,y)<\rho\}$. 

Towards a contradiction, suppose $U\cap T(U)$ is nonempty. Then there exists $u\in U$ for which $Tu\in U$. Then $d(Tu,x)<\rho\leq \frac{\epsilon}{2}$, and since $d(u,x)<\rho\leq \delta$, we also have $d(Tx,Tu)<\frac{\epsilon}{2}$. But then
\begin{align*}
    d(Tx,x)\leq d(Tx,Tu) + d(Tu,x) <\frac{\epsilon}{2}+\frac{\epsilon}{2}=\epsilon,
\end{align*}
which is a contradiction. This finishes the proof of the fact.
\end{proof}

There are many well-known examples of injective mixing operators on topological vector spaces. The most famous of these is due to Birkhoff in 1929, which in fact was the very first example of a hypercyclic operator \cite{Birkhoff2}.  Birkhoff showed that on the space of entire functions $H(\C)$, endowed with the usual topology of uniform convergence on compact subsets, any translation operator $T:H(\C)\to H(\C)$ defined by $T(f)(z)=f(z+a)$  is mixing, provided $a\in \C\setminus \{0\}$ (see \cite[Example 3.8]{GrossePeris}).  For our purposes, such examples show, along with Theorem \ref{injective}, that there are mixing operators which are not supermixing. 

We now turn our attention to strong topological transitivity. The definition of strong topological transitivity involves the zero vector, and thus it needs a slight modification to make sense for a general metric space $X$. We chose to strengthen the overall condition to the following: we say $T:X\to X$ is \textit{strongly topologically transitive} if, for each nonempty open subset $U$ of $X$, we have $X=\bigcup_{n=0}^\infty T^n(U)$. As stronger versions of the topological transitivity property, strong topological transitivity and weak-mixing are worth being compared to each other. It was shown in \cite[Example 4.3]{Ansari4} that a mixing mapping, which is necessarily weakly mixing, need not be strongly topologically transitive. Our next result shows that strong topological transitivity need not imply the weak-mixing property for general metric spaces.

\begin{theorem}\label{isometry theorem}
    Let $X$ be a complete metric space without isolated points and suppose $T:X\to X$ is a surjective, topologically transitive isometry. Then $T$ is strongly topologically transitive but not weakly mixing.
\end{theorem}

\begin{proof}
    By the Birkhoff transitivity theorem, which first appeared in \cite{Birkhoff1}, our assumptions on $X$ and $T$ imply $X$ must have a dense set of hypercyclic elements. Let $U\subseteq X$ be nonempty and open. Then there exists some element $a\in U$ with dense orbit and $\delta>0$ for which the ball $B_\delta(a)\subseteq U$. 

    Let $x\in X$. Since $\{a,Ta,T^2 a,\ldots\}$ is dense in $X$, there exists some $j\in \N$ such that $x\in B_\delta(T^j a)$. We show the image of the ball $B_\delta(a)$ under $T$ is given by $T(B_\delta(a))=B_\delta(Ta)$.
    
    First, let $y\in T(B_\delta(a))$. Then $Tx=y$ for some $x\in B_\delta(a)$, and since $T$ is an isometry, we have 
    \begin{align*}
        d(y,Ta)=d(Tx,Ta)=d(x,a)<\delta.
    \end{align*}
   Thus $y\in B_\delta(Ta)$, which shows $T(B_\delta(a))\subseteq B_\delta(Ta)$.

   Next, let $y\in B_\delta(Ta)$. Since $T$ is surjective, there exists $x\in X$ for which $Tx=y$. Then
    \begin{align*}
        d(x,a)=d(Tx,Ta)=d(y,Ta)<\delta,
    \end{align*}
    which shows $B_\delta(Ta)\subseteq T(B_\delta(a))$. Thus $T(B_\delta(a))=B_\delta(Ta)$.
    
    Since $T^j$ is also a surjective isometry, $x\in B_\delta(T^j a)$ implies $x\in T^j(B_\delta(a))$. Thus $x\in \bigcup_{n=1}^\infty T^n (B_\delta(a))\subseteq \bigcup_{n=1}^\infty T^n (U)$, which shows $T$ is strongly topologically transitive.

    We show $T$ is not weakly mixing by providing an open set $U$ for which $N(U,U)$ cannot contain two consecutive integers. Let $a\in X$ be a hypercyclic element for $T$, and let $\delta=d(Ta,a)$. Let $U=B_{\delta/4}(a)$. By means of contradiction, suppose $N(U,U)$ contains the integers $n$ and $n+1$. Then there exist $u_1,u_2\in U$ for which $T^n u_1\in U$ and $T^{n+1}u_2\in U$. Then
    \begin{align*} 
        d(Tu_2,u_1)=d(T^{n+1}u_2,T^n u_1)\leq d(T^{n+1}u_2, a)+d(a,T^n u_1)<\frac{\delta}{4}+\frac{\delta}{4}=\frac{\delta}{2},
    \end{align*}
    where we have used the fact that $T^n$ is an isometry to obtain the first equality above. But then
    \begin{align*}    
        d(Ta,a) &\leq d(Ta,Tu_2) + d(Tu_2,u_1) + d(u_1,a)\\
                &=d(a,u_2)+ d(Tu_2,u_1) + d(u_1,a)\\
                &< \frac{\delta}{4}+\frac{\delta}{2}+\frac{\delta}{4}=\delta,
    \end{align*}
    which is impossible.
\end{proof}

The quintessential application of the preceding theorem is to any ``irrational rotation" on the unit circle $S^1$ defined by $T:S^1\to S^1$, $Tz=e^{i 2\pi\theta} z$, where $\theta\in [0,1[$ is irrational. Such a circle rotation is easily seen to be a surjective isometry. That $T$ has a dense orbit takes more work to show; see \cite[Theorem 3.13]{Devaney} for a proof.

\section{Weighted backward shifts}\label{backward shifts}

The family of weighted backward shifts on $\ell^p$ and $c_0$ are the testing ground for any new notion in linear dynamics due to its simple yet rich structure. Recall that $c_0$ is the space of all bounded sequences $(x_n)_{n\geq 0}$ in $\C$ for which $x_n\to 0$, and the norm in $c_0$ is given by $\|x\|=\sup_{n\geq 0} |x_n|$. For $1\leq p<\infty$, $\ell^p$ is the space of all complex sequences $(x_n)_{n\geq 0}$ for which $\sum_{n=0}^\infty |x_n|^p <\infty$, and the norm on $\ell^p$ is given by $\|x\|=(\sum_{n=0}^\infty |x_n|^p)^{1/p}$.

Let $X$ be $c_0$ or one of the $\ell^p$ spaces, $1\leq p<\infty$. We denote by $e_n$, $n\geq 0$, the canonical basis vector of $X$ whose only nonzero term is a 1 in the $n$th position. Every $x=(x_0,x_1,\ldots)$ in $X$ can thus be written as $x=\sum_{i=0}^\infty x_i e_i$. If $w=(w_n)_{n\geq 0}$ is a bounded sequence of positive numbers, the weighted backward shift $B_w$ is defined by $B_w e_0=0$, $B_w e_n = w_n e_{n-1}$ for $n\geq 1$. We only consider surjective weighted backward shifts. It is a straightforward exercise to show that $B_w$ is surjective if and only if $\sup_{n\geq 0} \frac{1}{w_n} <\infty$. Each weighted backward shift has an associated forward shift $S$ defined on the basis of $X$ by $S e_n = \frac{e_{n+1}}{w_{n+1}}$. When $B_w$ is surjective, $S$ is defined on all of $X$, and the forward shift is a right inverse for $B_w$ on $X$, meaning $B_w S=I$, where $I$ is the identity operator on $X$. 

For simplicity in notations, we put $M_i^j = w_i w_{i+1}\cdots w_{i+j-1}$ for $i,j\geq 1$. Then for any $x\in X$, a formula for $S^n x$ is given by
\begin{align}\label{forward shift}
    S^n x = (0,0,\ldots,0,\frac{x_0}{M_1^{n}},\frac{x_1}{M_2^{n}}, \frac{x_2}{M_3^{n}},\ldots).
\end{align}
The respective characterizing conditions for $B_w$ to be topologically transitive (equivalently hypercyclic) or mixing are 
\begin{align*}
    \text{topologically transitive: } \sup_{n\geq 1} M_1^n = \infty, \hskip .25in \text{mixing: } \lim_{n\to \infty} M_1^n =\infty.
\end{align*}
For this result and much more background on backward shifts, we refer the interested reader to \cite[Section 4.1]{GrossePeris} and the references therein.

The characterizing condition on a weight sequence for $B_w$ to be hypermixing, which was obtained in \cite[Theorem 4.4]{Ansari4}, is
\begin{align}\label{hypermixing characterization}
    \sup_{n\geq 1} M_1^n = \infty \text{ and } \inf_{n,k\geq 1} M_n^k >0.
\end{align}
Observe the first condition above simply says $B_w$ is topologically transitive. The hypermixing characterization was obtained through the use of the following result  \cite[Theorem 4.1]{Ansari4}.
\begin{theorem}[{\cite[Theorem 4.1]{Ansari4}}]\label{Ansari theorem}
    For a weighted backward shift $B_w$ on $c_0$ or $\ell^p$, the following hold:
    \begin{enumerate}[label=\upshape(\roman*)]
        \item $B_w$ is hypermixing if and only if $B_w$ is surjective and $S^n x\to 0$ for every nonzero vector $x\in X$.
        \item $B_w$ is strongly topologically transitive if and only if $B_w$ is surjective and, for any nonzero vector $x\in X$, there is a sequence $(n_k)_k$ in $\N$ such that $S^{n_k} x\to 0$.
    \end{enumerate}
\end{theorem}

As mentioned above, a characterizing condition relying on the weight sequence for $B_w$ to be hypermixing is known, but such a characterization relying on the weight sequence for strong topological transitivity is not yet known. To discover this characterization, one would hope to generate a host of examples of strongly topologically transitive weighted backward shifts. Of course, any hypermixing $B_w$ is both strongly topologically transitive and mixing, and this gives us many examples already. Our Theorem \ref{shift theorem} shows that there are strongly topologically transitive weighted backward shifts which are not mixing. The proof of this result, combined with the results in \cite{Ansari4}, show the following diagram is complete.
\begin{diagram}
For the family of weighted backward shifts on $\ell^p$ or $c_0$, the following implications hold, and no implication that is not present holds in general.
    
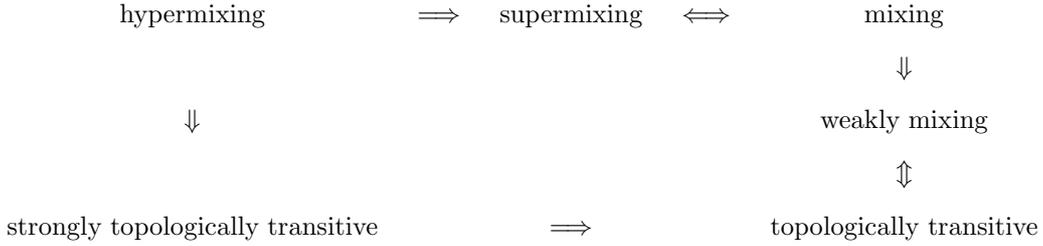
\begin{center}
\begin{tabular}{ccccc}
    hypermixing & $\implies$ & supermixing & $\iff$ & mixing\\
                &            &             &            & $\Downarrow$\\
    $\Downarrow$ &           &             &            & weakly mixing\\
                &            &             &            & $\Updownarrow$\\
    strongly topologically transitive &  & $\implies$ && topologically transitive
\end{tabular}
\end{center}
\end{diagram}

We are ready for the proof of our main result.
\begin{proof}[Proof of Theorem \ref{shift theorem}]
We construct a weight sequence $w=(w_n)_{n\geq 0}$ in a recursive fashion using blocks $b_1,b_2,b_3,$ etc. The even-numbered blocks will contain only $2$'s, and the odd-numbered blocks will contain only $\f$'s and $1$'s. Our first three blocks are defined as $b_1\coloneqq\{\f\}, b_2\coloneqq \{2\}$, and $b_3\coloneqq\{\f\}$. To describe the structure of the rest of our blocks, let $|b_n|$ represent the number of terms in the block $b_n$, and let $s_n=\sum_{i=2}^n |b_i|$, with $n\geq 3$. We define the rest of the blocks as follows. For even $n$, $b_n$ will be a block containing $s_{n-1}$ $2$'s. For odd $n \geq 5$, $b_n$ has the following properties:
\begin{enumerate}[label=\upshape(\roman*)]
    \item   $b_n$ contains the same number of $\f$'s as the number of $2$'s in the previous block, that number being exactly $|b_{n-1}|={s_{n-2}}$;
    \item $b_n$ begins and ends with $\f$;
    \item there are $s_{n-1}$ $1$'s between each successive pair of $\f$'s in $b_n$.
\end{enumerate}

Since the number of terms in each block increases very rapidly, we introduce some notation to help visualize the weight sequence.
Let $t^{(j)}$ represent the weight $t$ repeated $j$ times, and let $[t^{(j)} r]_k$ represent the string of weights $t^{(j)} r$ repeated $k$ times. For example, 
\begin{align*}
    \bigg[1^{(5)} \f\bigg]_3=11111\f11111\f11111\f.
\end{align*}
With this notation in hand, we can visualize our weight sequence.
\begin{align*}
    \underbrace{\f}_{b_1} \underbrace{ 2}_{b_2} \underbrace{\f }_{b_3} \underbrace{ 2^{(s_3)}}_{b_4} \underbrace{\f 1^{(s_4)} \f }_{b_5}  \underbrace{  2^{(s_5)}}_{b_6} \underbrace{\f \bigg[1^{(s_6)} \f\bigg]_{s_5-1}}_{b_7} \underbrace{ 2^{(s_7)}}_{b_8} \underbrace{\f \bigg[1^{(s_8)} \f \bigg]_{s_7 -1}}_{b_9} \cdots
\end{align*}

By property (i) above, multiplying together all weights in the blocks $b_2,b_3,\ldots b_{2k+1}$ will always yield 1. More specifically, since the combined number of terms in those blocks is $s_{2k+1}$, we have arranged the weight sequence so that $M_1^{s_{2k+1}}=1$ for all $k\geq 1$. Thus $B_w$ is not mixing since $M_1^n$ does not diverge to $\infty$. It remains to show that $B_w$ is strongly topologically transitive.

To show this, let $n_k=s_{2k+2}$ for each $k\in\N$. We claim the following inequalities hold for $k\geq 2$:
\begin{align}\label{easy est}
    M_i^{n_k} \geq \f \text{ whenever } i>s_{2k+1}, \text{ and }\\ \label{hard est}
    M_i^{n_k} \geq 2^k \text{ whenever } 1\leq i\leq s_{2k+1}.
\end{align}
Recall that $M_i^{n_k}$ is a product of $n_k$ successive weights beginning with the weight $w_i$. First, suppose $i>s_{2k+1}$. Since $s_{2k+1}$ is the total number of terms in blocks $b_2$ through $b_{2k+1}$, any weight $w_j$ with $j\geq i$ belongs to a block $b_l$ with $l\geq 2k+2$. By construction, $b_l$ either contains no weight equal to $\f$, or it contains $s_{2k+2}$ $1$'s separating any pair of $\f$'s that it contains. Thus the product $M_i^{n_k}$ contains at most a single factor of $\f$, and since the rest of the factors are either $1$ or $2$, (\ref{easy est}) follows.

Now suppose $1\leq i\leq s_{2k+1}$. We first underestimate the number of $2$'s in $M_i^{n_k}$. The first weight $w_i$ in the product $M_i^{n_k}$ must belong to one of the blocks $b_2,\ldots,b_{2k+1}$. Since $n_k$ is the total number of terms in blocks $b_2,\ldots,b_{2k+2}$, the product $M_i^{n_k}$ must contain each term from the block $b_{2k+2}$. Since $b_{2k+2}$ is a block consisting entirely of $2$'s, the product $M_i^{n_k}$ contains at least $|b_{2k+2}|=s_{2k+1}$ $2$'s.

We now overestimate the number of $\f$'s in $M_i^{n_k}$. Since the block $b_{2k+3}$ has $s_{2k+2}=n_k$ $1$'s between any successive pair of $\f$'s, the final weight in the product $M_i^{n_k}$ either appears in block $b_{2k+2}$ or in the first $n_k$ terms of block $b_{2k+3}$. By properties (ii) and (iii) of our construction, there are no $\f$'s in $b_{2k+2}$ and exactly one $\f$ in the first $n_k$ terms of block $b_{2k+3}$. Considering, by property (i), that the total number of $\f$'s in blocks $b_2$ through $b_{2k+1}$ is 
\begin{align*}
    1+|s_3|+\cdots+|s_{2k-1}|=|b_2|+|b_4|+\cdots+|b_{2k}|,
\end{align*}
an overestimate for the number of $\f$'s that appear in the product $M_i^{n_k}$ is
\begin{align*}
    1+|b_2|+|b_4|+\cdots+|b_{2k}|<s_{2k}.
\end{align*}
 Combining our underestimate of the number of $2$'s with our overestimate of the number of $\f$'s yields that $M_i^{n_k}$ is a power of $2$ whose power is at least $s_{2k+1}-s_{2k}=|b_{2k+1}|$. Since $|b_{2k+1}|\geq k$ when $k\geq 2$, the estimate (\ref{hard est}) follows.

We end by showing condition (ii) of Theorem \ref{Ansari theorem} is satisfied. Observe that $B_w$ is surjective since $\sup_{n\geq 0} \frac{1}{w_n}\leq 2$. Let $x=\sum_{i=0}^\infty x_ie_i$ be a nonzero vector in $\ell^p$ or $c_0$ and $\epsilon>0$. There exists some $N\in\N$ such that $\|\sum_{i=N}^\infty x_i e_i\|<\frac{\epsilon}{4}$. Choose $K$ such that for any $k\geq K$, we have both $\frac{1}{2^k}<\frac{\epsilon}{2\|x\|}$ and $s_{2k+1}> N$. Let $k> K$. If $i\geq s_{2k+1}$, then $M_{i+1}^{n_k}\geq \f$ by (\ref{easy est}), hence (\ref{forward shift}) implies 
\begin{align}\label{tail est}
    \bigg \|S^{n_k} \bigg(\sum_{i=s_{2k+1}}^\infty x_i e_i \bigg)\bigg\|\leq \bigg \| \sum_{i=s_{2k+1}}^\infty 2x_i e_i \bigg\| \leq \bigg \| \sum_{i=N}^\infty 2x_i e_i \bigg \|< \frac{\epsilon}{2}.
\end{align}
If instead $0\leq i<s_{2k+1}$, then $M_{i+1}^{n_k}\geq 2^k > \frac{2\|x\|}{\epsilon}$ by (\ref{hard est}), hence  (\ref{forward shift}) also implies 
\begin{align}\label{front est}
    \bigg\|S^{n_k} \bigg(\sum_{i=0}^{s_{2k+1}-1} x_i e_i\bigg)\bigg\|\leq  \bigg\| \sum_{i=0}^{s_{2k+1}-1} \frac{x_i}{2^k} e_i\bigg\| < \frac{1}{2^k} \|x\| <\frac{\epsilon}{2},
\end{align}
Together, estimates (\ref{tail est}) and (\ref{front est}) yield $\|S^{n_k} x\|< \epsilon$, which verifies condition (ii) of Theorem \ref{Ansari theorem}, completing the proof.
\end{proof}

The proof actually showed the existence of a sequence $n_k$ such that $S^{n_k}x\to 0$ for all $x\in X$, which is stronger than condition (ii) of Theorem 3.1 for strong topological transitivity. The proof of Theorem \ref{shift theorem} provides a new example of a strongly topologically transitive weighted backward shift. We invite the interested reader to pursue the following question.
\begin{question}
    Can we characterize the weight sequences which produce strongly topologically transitive weighted backward shifts on $\ell^p$ and $c_0$, like in (\ref{hypermixing characterization}) (which is \cite[Theorem 4.4]{Ansari4}) above?
\end{question}

Another natural question to ask was pointed out by an anonymous referee.
\begin{question}\label{big question}
    Does there exist a linear, bounded operator acting on a Banach space, which is strongly topologically transitive but not weakly mixing? 
\end{question}
\noindent That there exists a linear, bounded operator on a Banach space which is topologically transitive but not weakly mixing is a deep result due to De la Rosa and Read \cite{RosaRead}; see also \cite{BayartMatheron} for related examples on $\ell^p$ and $c_0$. One difficulty to answering the question is that such examples that would answer the question affirmatively, if they exist, are necessarily not weighted backward shifts, since every hypercyclic weighted backward shift satisfies the Hypercyclicity Criterion and is therefore weakly mixing (see \cite[Section 4.1]{GrossePeris}). 
The idea of the proof of Theorem \ref{isometry theorem} does not help, either, since no linear, bounded operator which is an isometry can be hypercyclic. Indeed, the orbit of any isometry is bounded, which precludes such an orbit from being dense in a Banach space.

\section{Acknowledgements}
We thank the anonymous referee for their corrections and suggestions which improved the quality of the paper. We also thank them for posing Question \ref{big question}.
\bibliography{Bib}

\end{document}